\documentclass{amsart}
\usepackage{amsmath}
\usepackage{amsfonts}
\usepackage{amssymb}

\newtheorem{theo}{Theorem}
\newtheorem{lem}[theo]{Lemma}
\newtheorem{prop}[theo]{Proposition}
\newtheorem{cor}[theo]{Corollary}
\newtheorem{defn}{Definition}

\newcommand{\E}{\ensuremath{\mathbb {E}}}
\renewcommand{\P}{\ensuremath{\mathbb {P}}}

\newcommand{\N}{{\mathbb N}}

\newcommand{\F}{{\mathcal F}}

\newcommand{\K}{{\mathcal K}}
\renewcommand{\H}{{\mathcal H}}

\newcommand{\C}{{\mathbb C}}

\def\E{\mathbb{E}}
\def\Z{\mathbb{Z}}

\def\F{{\mathcal F}}

\def\noi{\noindent}
\scrollmode

\subjclass[2010]{47A35, 47D37}

\scrollmode

\title[On N\"orlund summation and Ergodic Theory]{On N\"orlund summation and Ergodic Theory, with applications to 
power series of Hilbert contractions.}

\begin{document}

 \author{Christophe Cuny}
 
\address{Universit\'e de la Nouvelle-Cal\'edonie
Equipe ERIM,
B.P. 4477,
F-98847 Noumea Cedex}.

\email{cuny@univ-nc.nc}

\author{Michel Weber}

\address{IRMA, 10 rue du G\'en\'eral Zimmer,
67084 Strasbourg Cedex, France}
\email{michel.weber@math.unistra.fr}
\bigskip

\maketitle

\begin{abstract}  We show that if  ${\bf a}=(a_n)_{n\in \N}$ is a good weight for the dominated weighted ergodic theorem in 
$L^p$, $p>1$, then the N\"orlund matrix $N_{\bf a}=\{a_{i-j}/A_i\}_{0\le j\le i}$, $A_i=\sum_{k=0}^i |a_k|$ is bounded on $\ell^p(\N)$. 
We study the regularity (convergence in norm, almost everywhere) of operators in ergodic theory: power series of Hilbert contractions, and 
 power series $\sum_{n\in \N } a_nP_nf $ of $L^2$-contractions, and establish similar tight relations with  the N\"orlund operator associated to the modulus coefficient sequence $(|a_n|)_{n\in \N}$.
 \end{abstract}


\section{\bf Introduction}

Let ${\bf a}:= (a_n)_{n\in \N}$ be a  sequence of complex or real numbers 
(we take the convention $\N=\{0,1,2,\ldots\}$. 
We associate with 
${\bf a}$ an infinite matrix $N_{\bf a}=(a_{ij})_{i,j\in \N}$, called 
a N\"orlund matrix, in the following way. For every $i\ge 0$, set 
$A_i:=\sum_{k=0}^i |a_k|$.  We then define
\begin{equation}\label{defnorlund}
\begin{cases}a_{ij}:= a_{i-j}/A_i & \qquad \mbox{if $0\le j\le i$ and $A_i>0$,} \\
a_{ij}:= 0   &  \qquad \mbox{if $j>i$ or $A_i=0$}.\end{cases}\end{equation}


\noi Some authors consider instead $A_i=\sum_{k=0}^i a_k$, assuming then that it does not vanish.

\vskip 5pt

Then $N_{\bf a}$ induces naturally a (possibly unbounded) operator 
on $\ell^p(\N)$ for any $p\ge 1$. The matter of deciding whether this operator 
is bounded on some (or any)  $\ell^p(\N)$ is far from being 
solved. 
As noted by Bennett \cite{Bennett}, it seems, so far, that the best general known condition 
guaranteeing that $N_{\bf a}$ is bounded on any $\ell^{p}(\N)$, $p>1$, 
is that $a_n=O(A_n/n)$, see for instance Borwein and Cass \cite{BC}. That condition is realized when, for instance $(a_n)_{n\in \N}$ is a non-increasing sequence of positive numbers. 

\medskip

Let $c_{00}:=\{(u_n)_{n\in \N}\in \C^\N\, :\, \exists n_0\in \N\,  /\, 
 u_n=0,
\, \forall n\ge n_0\}$. Recall that the boundedness of $N_{\bf a}$ on $\ell^p$, $p> 1$, means that there 
exists $C_p>0$ such that, 
for any sequence $(u_n)_{n\in \N}\in c_{00}$,

\begin{equation}\label{boundedness}
\sum_{i\ge 0}|\frac1{A_i}\sum_{j=0}^ia_{i-j}u_j|^p \le C_p ^p
\sum_{i\ge 0}|u_i|^p\, .
\end{equation}

Equivalently, we have the dual formulation: for any sequence $(v_n)_{n\in \N}
\in c_{00}$, $q=p/(p-1)$,

\begin{equation}\label{boundedness2}
\sum_{j\ge 0}|\sum_{i\ge j}a_{i-j}v_i/A_i|^q\le C_p^q \sum_{j\ge 0}|v_j|^q\, .
\end{equation}
The latter is easily seen to be also equivalent to: for any sequence $(v_n)_{n\in \N}
\in c_{00}$,
\begin{equation}\label{boundedness3}
\sum_{j\ge 0}|\sum_{i\ge j}a_{i-j}v_i|^q
=\sum_{j\ge 0}|\sum_{i\ge 0}a_iv_{i+j}|^q\le C_p^q \sum_{j\ge 0}|A_jv_j|^q\, .
\end{equation} 
Moreover, it follows from \eqref{boundedness2} that
\begin{equation}\label{boundedness4}
\sum_{i\ge 0} \frac{|a_i|^p}{A_i^p}\le C_p^p\, ,
\end{equation}
where $a_i/A_i$ has to be interpreted as $0$ when $A_i=0$.

\medskip

We  show that N\"orlund matrices are connected with two different topics from ergodic theory. We establish tight relations between regularity (convergence in norm, almost everywhere) of operators in ergodic theory (power series of Hilbert contractions, power series of $L^2$-contractions, dominated weighted ergodic theorems, and naturally associated N\"orlund matrices. We   obtain conditions ensuring norm convergence of power series of Hilbert contractions, and also almost everywhere convergence of power series $\sum_{n\in \N } a_nP_nf $ of $L^2$-contractions. These conditions are expressed in terms of the N\"orlund operator associated to the modulus coefficient sequence $(|a_n|)_{n\in \N}$.

\section{\bf N\" orlund matrices and dominated weighted ergodic theorems} 

We first observe a connection between N\"orlund matrices 
and dominated weighted ergodic theorems. 

\medskip

We say that a sequence $(a_n)_{n\in \N}$, of complex numbers, is good for the 
dominated weighted ergodic theorem in $L^p$, $p>1$, if there exists $C>0$ such that 
for every 
dynamical system $(X,\Sigma,\nu,\tau)$, writing $A_n:=
\sum_{k=0}^{n}|a_k|$, we have
\begin{equation}\label{dom}
\|\sup_{n\ge 0} \frac1{A_{n}}|\sum_{k=0}^{n}a_kf\circ \tau^k|\, 
\|_{L^p(\nu)}\le C\|f\|_{L^p(\nu)} \qquad \forall f\in L^p(\nu)\,.
\end{equation}
Here again we take the convention that $\frac1{A_{n}}|\sum_{k=0}^{n}a_kf\circ \tau^k|=0$ if 
$A_n=0$.

 The next lemma is well-known, it is in the spirit of the so-called Conze principle, see for instance \cite[Th.\,5.4.3]{Weber}. It states a converse of Calderon's transference principle. 

\begin{lem}
Let $(a_n)_{n\in \N}$ be good for the dominated weighted ergodic theorem in 
\marginpar{$p>1$}$L^p$, $p>1$. Then, with the best constant $C>0$ appearing in \eqref{dom}, we have
for every $(v_n)_{n\in \Z}\in \ell^p(\Z)$, 
\begin{equation}\label{dominatedshift}
\sum_{i\in \Z} \Big( \sup_{n\ge 0} \frac1{A_{n}}|\sum_{j=0}^{n} a_jv_{i+j}
|\Big)^p \le C^p \sum_{i\in \Z} |v_j|^p\, .
\end{equation}
\end{lem}
\noindent {\bf Proof.} 
Let $(v_n)_{n\in \Z}\in \ell^p(\Z)$. Let $N>M\ge 1$ be integers (one has in mind that $N\gg M$). Take $X=\{-N,-N+1,\ldots, N-1,N\}$,  $\nu:=\frac1{2N+1}\sum_{k=-N}^N\delta_k$ and $\tau$ the transformation given by $\theta(k)=k+1$ if $k\neq N$ and $\theta(N)=-N$. Define $f$ on $X$, by $f(k)=v_k$ for every $k\in X$.
By \eqref{dom}, we have 
\begin{gather}
\nonumber\frac1{2N+1}\sum_{i=-N}^{N+1-M} \Big(\sup_{0\le m\le M}\frac{1}{A_{m}
} |\sum_{k =0}^{m}a_k v_{i+k}|\Big)^p 
\le  \|\sup_{n\ge 0} \frac1{A_{n}}|\sum_{k=0}^{n}a_kf\circ \tau^k|\, 
\|_{L^p(\nu)}^p\\
\label{maj}\le C^p\|f\|_{L^p(\nu)}^p=  \frac{C^p}{2N+1}\sum_{i=-N}^N |v_i|^p\, .
\end{gather}
Multiplying \eqref{maj} by $2N+1$, letting first $N\to +\infty$ and then $M\to +\infty$, we derive \eqref{dominatedshift}. \hfill $\square$ 

\noindent {\bf Remark.} Our proof is based on the use of the dominated 
weighted ergodic theorem on periodic systems (the rotations on $\Z/(2N+1)\Z$).  
To give a proof based on the dominated weighted ergodic theorem 
on a \emph{single} (but ergodic and non-atomic) dynamical system, one could 
use Rohlin's lemma (see for instance Weber \cite[p.\,270]{Weber} 
for a statement of the lemma). 

\medskip

We deduce the following.

\begin{prop}\label{prop-erg-norl}
Let ${\bf a}=(a_n)_{n\in \N}$ be a good weight for the dominated weighted ergodic theorem in 
 $L^p$, $p>1$. Then, the N\"orlund matrix $N_{\bf a}$ is bounded on $\ell^p(\N)$. 
Moreover, for every non-increasing sequence of nonnegative numbers 
$(b_n)_{n\in \N}$, writing ${\bf c}:=(a_nb_n)_{n\in \N}$, 
$N_{\bf c}$ is bounded on $\ell^p(\N)$.
\end{prop}
\noindent {\bf Remark.} It is unclear whether "$N_{\bf a}$ bounded on $\ell^p$" implies   "$N_{\bf c}$ bounded on $\ell^p$", in general. 

\noindent {\bf Proof.} 
Let $(u_n)_{n\in \N}\in \ell^p(\N)$. Define $(v_n)\in \ell^p(\Z)$ as follows. 
$v_n:=u_{-n}$ if $n\le 0$ and $v_n:=0$ if $n>0$. Using \eqref{dominatedshift} 
and for every $i\ge 1$ the trivial estimate 
$$
\frac1{A_{i} }|\sum_{j=0}^{i} a_jv_{-i+j}|\le \sup_{n\ge 0} 
\frac1{A_{n} }|\sum_{j=0}^{n} a_jv_{-i+j}|\, ,
$$
we infer that 
$$
\sum_{i\ge 0} \Big( \frac1{A_{i} }|\sum_{j=0}^{i} a_jv_{-i+j}|\Big)^p\le 
C^p \sum_{i\in \Z} |v_i|^p =C^p \sum_{i\ge 0}|u_i|^p\, .
$$
Using that $v_{-i+j}=u_{i-j}$ when $ j\le i$, we derive that 
$N_{\bf a}$ is bounded on $\ell^p$. 

To prove the last assertion, one just has to notice that, using 
 Abel summation, $(c_n)_{n\in \N}$ is a good weight for the dominated weighted ergodic theorem. \hfill $\square$. 

\medskip

Of course, as one can see from the above proof, the fact that ${\bf a}=(a_n)_{n\in \N}$ be a good weight for the dominated weighted ergodic theorem in $L^p$ is a much stronger statement than the fact that 
  $N_{\bf a}$ be bounded on $\ell^p(\N)$. Hence, Proposition 
  \ref{prop-erg-norl} should not be seen as a method to prove boundedness 
  of some N\"orlund matrices, but as a source of examples of 
  N\"orlund matrices, since there are many examples of sequences that are 
  known to be good for the dominated weighted ergodic theorem. 
  We provide some of them below. One may also consult the survey 
  by Bellow and Losert \cite{BL} for dominated weighted ergodic theorems 
  with bounded weights. More arithmetical sequences may be found 
  in Cuny and Weber \cite{CW}.
  
  \medskip
  
 \noindent {\bf Examples.} The following sequences $(a_n)_{n\in \N}$ are 
 good for the dominated weighted ergodic theorem in $L^p$, for every 
 $p>1$:

 $(i)$ (Bourgain and Wierdl, \cite{Bourgain-prime}, \cite{Wierdl}) Let 
 ${\mathcal P}$ be the set of prime numbers and take $a_n:=
 \delta_{\mathcal P}(n)$, for every $n\in \N$. 
 
 $(ii)$ (Bourgain, \cite{Bourgain-square}) Let 
  ${\mathcal S}$ be the set of squares and take $a_n:=
 \delta_{\mathcal S}(n)$, for every $n\in \N$.
 
 $(iii)$ (Cuny and Weber, \cite{CW}) Take $a_0=0$ and for every $n\in \N$ take $a_n=d_n$, the number of divisors of $n$. 
 \vskip 4 pt 
We now give an example which does not work on every $L^p$, $p>1$. 
Let $(X,\Sigma,\mu,\theta)$ be an ergodic dynamical system. 
Let $g\in L^q(\mu)$, for some $1< q\le \infty$. 

\smallskip
$(iv)$ (Bourgain, Demeter, Lacey, Tao and Thiele, \cite{Bourgain-return}, 
\cite{DLTT} and \cite{Demeter}) There exists 
$\overline X\in \Sigma$ with $\mu(\overline X)=1$ such that for every 
$x\in \overline X$, $(a_n)_{n\in \N}:=(g\circ \theta^n(x))_{n\in \N}$ 
is good for the dominated weighted ergodic theorem in $L^p$ for every 
$p>1$ such that $1/p+1/q<3/2$.

\medskip

Let us notice that none of the above examples satisfies  the previously mentionned criterium: $\sup_{n\in \N}n|a_n|/A_n<\infty$. The fact that the N\"orlund matrix associated with the sequence $(a_n)_{n\in \N}$ in example $(ii)$ is bounded  has been proved by Borwein 
\cite{Borwein}. 

\section{\bf Norm convergence of power series of Hilbert contractions}
\label{section-norm} 

Let $P$ be a contraction of a (real or complex) Hilbert space 
$\H$. Given a sequence of complex 
numbers $(a_n)_{n\in \N}$ and $f\in \H$,  we are interested in finding 
 conditions involving $(\|P^nf\|_\H)_{n\in \N}$ sufficient for the 
norm convergence of $\sum_{n\in \N}a_nP^nf$. 

\smallskip

An obvious condition is the following 

\begin{equation}\label{obvious}
\sum_{n\in \N}|a_n|\|P^nf\|_\H<\infty\, .
\end{equation}

Sufficient conditions involving $(\|f+\ldots +P^nf\|_\H)_{n\in \N}$ 
have been obtained when $P$ is unitary (i.e. $P^*=P^{-1}$) or, more generally, normal  (i.e. $PP^*=P^*P$), if moreover $(a_n)_{n\in \N}$ 
is regular (at least nonnegative and nonincreasing). Let us mention 
the papers \cite{Gaposhkin-series} and \cite{Cuny-series}, see also \cite{Cuny-Lp} for some 
$L^p$ versions.

Recall that, see for instance Nagy and Foias \cite{NF} (see also Sh\"affer \cite{shaffer} for an explicit matrix construction), $P$ admits a unitary dilation, that is, there exist another Hilbert space $\K$, with $\H\subset 
\K$, and a unitary operator $U$ on $\K$ such that $EU^n=P^n$ for every 
$n\ge 1$, where $E$ is the orthogonal projection onto 
$\H$.

\medskip 

We start with some simple lemmas. The first one appears in 
Cuny and Lin \cite{CL}, but we recall the short proof. 

\begin{lem}\label{ortho}
For every $n\in \N$ and every $\ell\ge 1$, the spaces 
$(U^{-n}P^n -U^{-n-1}P^{n+1})\H$ and $U^{-n-\ell}P^{n+\ell}\H$ are 
orthogonal (in $\K$).
\end{lem}
\noindent {\bf Proof.} Let $f,g\in  \H$. Let $n\in \N$ and $\ell\ge 1$. We have 
\begin{gather*}
\langle (U^{-n}P^n -U^{-n-1}P^{n+1}) f, U^{-n-\ell}P^{n+\ell} g\rangle_\K
=\langle U^{\ell} P^nf ,P^{n+\ell} g\rangle_\K -\langle U^{\ell-1} 
P^{n+1}f ,P^{n+\ell} g\rangle_\K\\ 
=\langle P^{n+\ell} f ,P^{n+\ell} g\rangle_\K-\langle P^{n+\ell} f ,P^{n+\ell} g\rangle_\K=0\, .
\end{gather*}
\hfill $\square$

\begin{lem}\label{lemeq}
Let $f\in \H$ be such that $\|P^mf\|_\H \to 0$ as $m\to +\infty$. 
Then, for every $n\ge 1$,  $\|P^nf\|_\H^2 =\|P^nf\|_\K^2  = \sum_{k\ge n} \|U^{-k}P^{k}f -U^{-k-1}P^{k+1}\|_\K^2$. In particular, for any positive and non-decreasing sequence $(b_n)_{n\ge 0}$, the following  are 
equivalent (setting $b_{-1}=0$). 
\begin{itemize}
\item [$(i)$] $\sum_{n\in \N} (b_{n}-b_{n-1})\|P^nf\|_\H^2 <\infty$;
\item [$(ii)$]  $\sum_{n\in \N} b_{n}\|U^{-n}P^nf- 
U^{-n-1}P^{n+1}f\|_\K^2 <\infty$.
\end{itemize}
\end{lem}
\noindent {\bf Remarks.} Notice that by Kronecker's lemma, if $(i)$ holds 
$\|P^nf\|_\H^2\sum_{k=0}^n(b_{k+1}-b_k)=\|P^nf\|_\H^2(b_{n+1}-b_0)
\to 0$ as $n\to +\infty$. In particular, since $(b_n)_{n\in \N}$ is non decreasing, $\|P^nf\|_\H\to 0$. Item $(i)$ is satisfied if $\sum_{n\ge
 0}b_{2^{n+1}}\|P^{2^n}f\|_\H^2 <\infty$.
\medskip

\noindent {\bf Proof.} Since $\|P^nf\|\to 0$,  for every 
$n\in \N$, we have, (with convergence in $\K$)
\begin{gather}\label{Pnf}
P^n f =\sum_{k\ge 0} (U^{-k}P^{n+k}f-U^{-k-1}P^{n+k+1}f)\, .
\end{gather}
By the above lemma the terms of that series lie in orthogonal spaces. 
Hence, 
\begin{align*}
\|P^nf\|_\K^2  &  = \sum_{k\ge 0} \|U^{-k}P^{n+k}f -U^{-k-1}P^{n+k+1}\|_\K^2 \\
               &  = \sum_{k\ge n} \|U^{-k}P^{k}f -U^{-k-1}P^{k+1}\|_\K^2 \, ,
\end{align*}
where we used that $U$ is unitary (and a change of variable) 
for the last identity. Then, the equivalence of  $(i)$ and $(ii)$ follows by Fubini. \hfill $\square$


\smallskip

Given a sequence of complex numbers $(a_n)_{n\in \N}$, consider the following conditions 

\begin{gather}
\label{suff-cond-1} \sum_{n\in \N}|a_n|\big(\sum_{k=0}^n |a_k|
\big)\|P^nf\|_\H 
^2 <\infty\, ,\\
\label{suff-cond-2}\sum_{n\in \N} \big(\sum_{k=0}^n |a_k|\big)^2\|U^{-n}P^nf- 
U^{-n-1}P^{n+1}f\|_\K^2 <\infty
\end{gather}

By Lemma \eqref{lemeq}, when $\|P^nf\|_\H\to 0$, \eqref{suff-cond-1} and \eqref{suff-cond-2} 
are equivalent. Assume that \eqref{obvious} holds. Then, since 
$(\|P^nf\|_\H)_{n\in \N}$ is nonincreasing, $\sup_{n\in \N} 
\|P^nf\|\sum_{k=0}^n |a_k| <\infty$ and \eqref{suff-cond-1} holds. Hence, \eqref{suff-cond-1} is always weaker than \eqref{obvious}.

\begin{prop}\label{prop-norm}
Let $(a_n)_{n\in \N}\in \C^\N$ be such that $N_{\bf |a|}$ be bounded 
on $\ell^2( \N)$ where ${\bf |a|}=(a_n)_{n\in \N}$. Let $f\in \H$ be such 
 that either of conditions 
\eqref{suff-cond-1} or \eqref{suff-cond-2} hold. Then, the series
$\sum_{n\in \N}a_nP^n f$ converges in $\H$.
\end{prop}
\noindent {\bf Proof.} 
Since  $N_{\bf |a|}$ is  bounded 
on $\ell^2( \N)$, then by \eqref{boundedness4} (with $p=2$)
\begin{equation}\label{remnor}
\sum_{n\in \N} \frac{a_n^2}{A_n^2}<\infty\, .
\end{equation}
Let $q>p\ge 1$ be integers and write 
$V_{p,q}f:=\sum_{k=p}^q a_k P^kf$. For every $n\in \N$, let $u_n:=\|U^{-n}P^{n}f-U^{-n-1}P^{n+1}f\|_\K$ and $v_n:=A_nu_n$, where $A_n=
\sum_{k=0}^n |a_k|$. Finally, let ${\bf v}:=(v_n)_{n\in \N}$. By Lemma \ref{lemeq} and using that 
$U$ is unitary, we have
\begin{gather}\label{increment}
\|V_{p,q}f\|_\K^2 =\sum_{n\in \N}\|U^{-n}P^nV_{p,q}f-
U^{-n-1}P^{n+1}V_{p,q}f\|_\K^2  \le 
\sum_{n\in  \N}\Big( \sum_{k=p}^q |a_k|u_{n+k}\Big)^2 \, .
\end{gather}
By Cauchy's criteria one has to prove that $\|V_{p,q}f\|_\K
\to 0$ as $p,q\to +\infty$. Using the Lebesgue dominated theorem for the counting measure on $\N$, it suffices to prove that 
\begin{equation}\label{cauchy}
\sum_{k=p}^q |a_k|u_{n+k}\underset{p,q\to +\infty}\longrightarrow 
0\, ,
\end{equation}
and that 
\begin{equation}\label{domination}
\sum_{n\in  \N}\Big( \sum_{k\ge 0} |a_k|u_{n+k}\Big)^2\, .
\end{equation}
The convergence  \eqref{cauchy}  follows from Cauchy-Schwarz combined with the assumed condition 
\eqref{suff-cond-2} and \eqref{remnor}. 

FTo prove \eqref{domination}, it suffices to notice that 
\begin{gather*}
\sum_{n\in  \N}\Big( \sum_{k\ge 0} |a_k|u_{n+k}\Big)^2 \le 
\sum_{n\in  \N}\Big( \sum_{k\ge n}|a_{k-n}|u_{k}\Big)^2
=\|N_{\bf |a|}^*{\bf v}\|_{\ell^2(\N)}^2\\ \le 
\|N_{\bf |a|}^*\|^2\, \|{\bf v}\|_{\ell^2(\N)}^2=
\|N_{\bf |a|}\|^2\, \sum_{n\in \N}A_n^2u_n^2\, .
\end{gather*}
\hfill $\square$

\medskip 
The proposition has been proved in \cite{CL} in the case where 
$a_n=n^{-1/2}$. An important case corresponds to the situation 
where $a_n=1$ for every $n\in \N$. Then, the proposition 
gives a sufficient condition (namely $\sum_{n\in \N}n\|P^nf\|_\H^2<\infty$) 
for $f$ to be a coboundary (i.e. $f=(I-P)g$ for some $g\in \H$). 
This sufficient condition has been obtained independently by Voln\'y \cite{Volny} in the special case 
where $P$ is a Markov operator on $L^2(m)$. His proof  (which 
does not appeal to the notion  of N\"orlund matrices) is essentially the 
same, since the shift on the space of trajectories of the associated Markov chain plays the role of the unitary dilation. 

\medskip

\begin{prop}\label{prop-eq}
Let $(a_n)_{n\in \N}\in \C^\N$. Assume that for every contraction $P$ 
on a Hilbert space $\H$ the following property holds : "If \eqref{suff-cond-1} holds for some $f\in \H$ then $\sum_{n\in \N}a_n P^nf$ 
converges in $\H$". Then, $N_{\bf a}$ is bounded on  $\ell^2$. 
\end{prop}
\noindent {\bf Proof.} Let $P$ be a contraction on a Hilbert space $\H$
satisfying the above property. Let ${\mathcal L}:=\{f\in \H\,:\, 
\sum_{n\in \N}|a_n|\big(\sum_{k=0}^n|a_k|\big)\|P^nf\|_\H^2<\infty$. Then, ${\mathcal L}$ is 
a Hilbert space and we define an operator $T$ on ${\mathcal L}$, by 
setting $Tf=\sum_{n\in \N}a_nP^nf$ for every $f\in {\mathcal L}$.   Then, by  the Banach-Steinhaus, theorem $T$ is continuous. Hence, there exists $C=C_{\H,P}$, such that 
$\|\sum_{n\in \N} a_n P^nf\|_\H\le C( \sum_{n\in \N}(b_{n+1}-b_n)
\|P^nf\|_\H^2)^{1/2}
<\infty$. 

\smallskip

Let us prove the proposition. We give a probabilistic proof. Let $(\Omega,\F,\P)$ be the  probability space given by $\Omega=\{-1,1\}^\Z$, $\F$ the product $\sigma$-algebra and $\P=\mu^{\otimes \Z}$, with $\mu(0)=\mu(1)=1/2$.  Let 
$\theta$ be the shift on $\Omega$ and $(\varepsilon_n)_{n\in \Z}$ 
be the coordinate process. In particular, $\varepsilon_{n+1}=\varepsilon_n 
\circ \theta$ and $(\varepsilon_n)_{n\in \Z}$ is iid. 

\medskip

Denote $\F_0:=\sigma\{\varepsilon_i,i\le 0\}$. 
Set $\H:=L^2(\Omega,\F_0,\P)$ and $\K:=L^2(\Omega,\F,\P)$ and define 
two operators $U$ and $P$ on $\K$ and $\H$ respectively by 
$Uf=f\circ \theta$ for every $f\in \K$ and $Pf=\E(f\circ \theta|\F_0)$ 
for every $f\in \H$ (then $P$ is a Markov operator). Clearly, $U$ is a unitary dilation of $P$. Let $(u_i)_{i\in \N}
\in c_{00}$ and define $f:=\sum_{i\in \N}u_i\varepsilon_i \in\H$. 
Assume moreover that $\sum_{n\in \N}|a_n|\big(\sum_{k=0}^n|a_k|\big))\|P^nf\|_\H^2
<\infty$, or equivalently (by Lemma \ref{lemeq}), 
$\sum_{n\in \N}\big(\sum_{k=0}^n|a_k|\big)^2 \|U^{-n}P^n f-U^{-n-1}P^{n+1}f\|_2^2$. 
Notice that $P^nf= \sum_{i\in \N} u_{i+n}\varepsilon_{-i}$ and that 
$ \|U^{-n}P^n f-U^{-n-1}P^{n+1}f\|_2^2 =u_{n}^2$. Moreover, 
$\|\sum_{n\in \N}a_nP^nf\|_\H^2 =\sum_{i\in \N} |\sum_{n\in \N}a_n
u_{i+n}|^2$. Hence, $\sum_{i\in \N} |\sum_{n\in \N}a_n u_{i+n}|^2
\le C^2 \sum_{n\in \N}b_n|u_n|^2$, i.e. \eqref{boundedness3} holds 
with $q=2$, and the proof is complete. 
\hfill $\square$

\medskip

We shall now prove that Proposition \ref{prop-norm} cannot be improved. 

\begin{defn}
We say that a contraction $P$ on $\H$ is Ritt if 
$\sup_{n\in \N} n\|P^n-P^{n+1}\|<\infty$.
\end{defn}

\begin{prop}\label{prop-alpha}
Let $P$ be a  contraction on $\H$. For every $0\le \alpha <1$, 
consider the following properties. 
\begin{itemize}
\item [$(i)$] The series $\sum_{n\in \N} (n+1)^{-\alpha}P^nf$ converges in $\H$ ;
\item [$(ii)$] $\sum_{n\in \N} (n+1)^{1-2\alpha}\|P^nf\|_\H^2<\infty$\, .
\end{itemize}  
Then, $(ii)\Rightarrow (i)$. If moreover $P$ is Ritt then 
$(i)\Rightarrow (ii)$. 
\end{prop}
\noindent {\bf Remark.} By \cite{CCL2}, when $P$ is a positive operator on $L^2(m)$ then $(i)$ of the proposition implies that the series 
$\sum_{n\in \N} (n+1)^{-\alpha}P^nf $ converges $m$-almost everywhere and the associated maximal function is in $L^2$. The fact that $(i)\Rightarrow (ii)$ 
has been proved  by Cohen, Cuny and Lin \cite{CCL} using results from Arhancet and Le Merdy \cite{AL}) when $\alpha\in (0,1)$ and 
$P$ is a positive Ritt contraction of some $L^2(m)$ (there are also 
 analogous results in $L^p$ in \cite{CCL}). 
 

\medskip

\noindent {\bf Proof.} The fact that $(ii)\Rightarrow (i)$ is a direct application of Proposition \ref{prop-norm}. Assume that $P$ is a Ritt operator and that $\sum_{n\in \N} \frac{P^nf}{(n+1)^\alpha}$ converges in $\H$.

We start with the case $0<\alpha<1$. By Proposition 4.6 of Cohen, Cuny and Lin \cite{CCL2} 
(see also their example (v) page 8), we have 
\begin{gather*}
\sum_{n\ge 0} \frac{\|Pf+\cdots +P^{2^n}f\|_\H^2}{2^{2\alpha n}}
<\infty \, ,
\end{gather*}
Then, using (3) of Cohen, Cuny and Lin \cite{CCL} 
combined with Lemma \ref{subadditive} below, we infer that 
$\sum_{n\ge 0}2^{(2-2\alpha)n}\|P^{2^n}f\|_\H^2 <\infty$, which finishes 
the proof, in that case. 

Assume now that $\alpha=0$. Let $g:=\sum_{n\in \N} 
P^nf$. Then, $f=(I-P)g$. Hence, by Theorem 8.1 of 
Le Merdy \cite{Lemerdy}, 
$$
\sum_{n\in \N} n\|P^nf\|_\H^2=\sum_{n\in \N} n\|P^n(I-P)g\|_\H^2
\le \|g\|_\H\, ,
$$ 
which is the desired result.  \hfill $\square$
\medskip

\section{\bf Almost everywhere convergence of power series of 
$L^2$-contractions}\label{section-ae}

Once norm convergence has been proven, one may wonder, in the case where $\H=L^2(m)$, 
 whether almost everywhere convergence holds. As mentionned in the remark following Proposition 
 \ref{prop-alpha}, for "regular" sequences, if $P$ is a positive contraction of $L^2(m)$ 
 then norm convergence implies almost everywhere convergence. However, as we shall see below 
 (see Proposition \ref{Cex}), there is no such  result for contractions that are not positive. Let us mention that 
 the almost everywhere convergence of power series (for regular $(a_n)_{n\in \N}$) for unitary or normal operators 
 on $L^2(m)$ has been proven under conditions involving $(\|f+\ldots +P^nf\|_\H)_{n\in \N}$ 
 in \cite{Gaposhkin-series} and \cite{Cuny-series}, see also \cite{Cuny-Lp2} for 
 $L^p$-versions. 

\begin{theo}\label{theo-ae}
Let $(a_n)_{n\in \N}\in \C^\N$ be such that $N_{\bf |a|}$ be bounded 
on $\ell^2( \N)$ where ${\bf |a|}=(|a_n|)_{n\in \N}$. Let $A_n:=\sum_{
k=0}^n |a_k|$. Let $P$ be a contraction on $L^2(m)$. 
Let $f\in L^2(m)$ such that
\begin{equation}\label{cond-ae}
\sum_{n\ge 1}  ( \log (n+1))^2 A_{2^{n+1}}^2
\|P^{2^n}f\|_{L^2(m)}^2<\infty\, .
\end{equation}
Then, the series $\sum_{n\in \N}a_n P^n f$ converges $m$-almost everywhere 
and 
$$\sup_{N\ge 1}\Big|\sum_{n=0}^N a_nP^nf\Big|\in {L^2(m)}.$$
\end{theo}
\noindent {\bf Remark.} A sufficient condition for \eqref{cond-ae} is the following
\begin{equation}\label{cond-ae-bis}
\sum_{n\ge 1}  ( \log\log (n+3))^2 \frac{A_{4n}^2}{n+1}
\|P^{n}f\|_{L^2(m)}^2<\infty\, .
\end{equation}

\noindent {\bf Proof.} 
Let $N\in \N$. We have 
\begin{gather*}
\max_{2^N\le n\le 2^{N+1}-1}|\sum_{k=2^N}^n a_n P^kf|\le 
\sum_{k=2^N}^{2^{N+1}-1}|a_k|\, |P^k f|\, .
\end{gather*}
Hence, 
\begin{gather*}
\sum_{N\in \N} \|\, \max_{2^N\le n\le 2^{N+1}-1}|\sum_{k=2^N}^n a_n P^kf|\, \|_{L^2(m)}^2
\le \sum_{N\in \N} A_{2^{N+1}}^2\|P^{2^N}f\|_{L^2(m)}^2<\infty\, .
\end{gather*}

In particular, it suffices to prove that $(\sum_{n=0}^{2^N}a_nP^nf)_{N\ge 0}
$ converges and that $\sup_{N\ge 0}|\sum_{n=0}^{2^N}a_nP^nf|
\in L^2(m)$. 

\medskip

By \eqref{increment}, for every $q\ge p$, we have 
\begin{equation}\label{increment-2}
\|\sum_{n=2^p}^{2^q-1}a_n P^nf\|_{L^2(m)}^2\le \sum_{n\in \N}\Big(\sum_{k=2^p}
^{2^q-1}|a_k|u_{n+k}\Big)^2\, .
\end{equation}
Set $d(p,q):= \sum_{n\in \N}\Big(\sum_{k=2^p}
^{2^q-1}|a_k|u_{n+k}\Big)^2$ and notice that $d$ is super-additive in the following sense: for every $m\ge l\ge k$,  $d(k,l)+d(l,m)\le d(k,m)$. 
By Proposition 2.2 of Cohen and Lin \cite{CoL}, there exists $C>0$, such that 
for every $n\ge 0$, 
$$
\|\max_{2^{2^n}\le m\le 2^{2^{n+1}}-1}|\sum_{k=2^{2^n}}^m a_k P^kf|\, 
\|_{L^2(m)}^2\le C (n+1)^2 d(2^n,2^{n+1}-1)\, .
$$
Assume that 
\begin{equation}\label{inter-cond}
\sum_{n\ge 0}(n+1)^2 d(2^n,2^{n+1}-1) <\infty\, .
\end{equation}
Then, using \eqref{increment-2} and Cauchy-Schwarz we see that
$$
\big(\sum_{n\in \N}\|\sum_{k=2^{2^n}}^{2^{2^{n+1}}-1}a_k P^k f\|_{L^2(m)}\big)^2
\le \sum_{n\in \N}\frac{1}{(n+1)^2}\sum_{n\in \N}(n+1)^2
d(2^n,2^{n+1}-1)\, .
$$
This finishes the proof, provided that we can show \eqref{inter-cond}. 

\medskip

But, \eqref{inter-cond} reads 
$$
\sum_{n\in \N}\sum_{\ell\ge 0}(\ell +1)^2\Big(\sum_{k=2^{2^\ell}}
^{2^{2^{\ell+1}}-1}|a_k|u_{n+k}\Big)^2<\infty \, .
$$
Using that $\|\cdot \|_{\ell^2}\le \|\cdot \|_{\ell^1}$, we infer that 
\begin{gather*}
\sum_{n\in \N}\sum_{\ell\ge 0}(\ell +1)^2\Big(\sum_{k=2^{2^\ell}}
^{2^{2^{\ell+1}}-1}|a_k|u_{n+k}\Big)^2\le 
\sum_{n\in \N}\Big(\sum_{k\ge 0}
(\log\log (k+3))^2|a_k|u_{n+k}\Big)^2\\
\le 
\sum_{n\in \N}\Big(\sum_{k\ge 0}
(\log\log (n+k+3))^2|a_k|u_{n+k}\Big)^2\, .
\end{gather*}
Then, proceeding as in the (end of the) proof of Proposition 
\ref{prop-norm} we see that \eqref{inter-cond} holds provided that 
$$
\sum_{n\in \N} (\log\log (n+3))^2 A_n^2 u_n^2<\infty\, ,
$$
which follows from \eqref{cond-ae} using that $(A_n)_{n\in \N}$ is non-decreasing and that $\sum_{k=2^n}^{2^{n+1}-1}u_k^2\le \|P^{2^n}f\|_{L^2(m)}
^2$. \hfill $\square$

\vskip 5 pt
\begin{cor}\label{cor}
Let $(X,\Sigma,\mu,\theta)$ be an ergodic dynamical system. 
Let $g\in L^p(\mu)$ for some $p>1$. There exists 
$\overline X\in \Sigma$ with $\mu(\overline X)=1$ such that for every 
$x\in \overline X$, setting $(a_n)_{n\in \N}:=(g\circ \theta^n(x))_{n\in \N}$ the following holds: for every  $0\le \alpha <1$, every contraction $P$ on  $L^2(m)$ and every $f\in L^2(m)$  such that 
$$
\sum_{n\in \N} (\log \log (n+3))^2(n+1)^{1-2\alpha}\|P^nf\|_2^2
<\infty\, ,
$$
the sequence   $\sum_{n\in \N}\frac{a_n P^nf}{(n+1)^\alpha}$ 
converges $m$-almost everywhere and the associated maximal function if in $L^2(m)$.  
\end{cor}
 
\noindent {\bf Proof.}
Let $(X,\Sigma,\mu,\theta)$ and let $g\in L^p(\mu)$. 
Let $\overline X$ be the set appearing in the example $(iv)$. 
Modifying $\overline X$ if necessary we may assume that 
$A_n=|a_0|+\ldots + |a_n|\le K(x) n$, for some finite $K(x)>0$. 
Then, for every $x\in \overline X$, $(g\circ \theta^n(x))_{n\in \N}$ 
is good for the dominated weighted ergodic theorem. Applying 
Proposition \ref{prop-erg-norl}, we see that, with 
${\bf c}=(c_n)_{n\in \N}:=((n+1)^{-\alpha}f\circ 
\theta^n(x))_{n\in \N}$, $N_{\bf c}$ is bounded on $\ell^2$. 
Set $C_n:= \sum_{k=0}^n |c_k|$ (we see $C_n$ as a function on $X$). By Theorem \ref{theo-ae} (see the remark after the theorem), we
are back to prove that $\sum_{n\ge 1}(\log\log (n+1))^2 
C_{4n}^2\|P^nf\|_2^2<\infty$. 
But this follows  our assumption (and an Abel summation) since $A_n\le K(x)n$. \hfill $\square$

\medskip

We shall now prove that Corollary \ref{cor} (and hence Theorem \ref{theo-ae}) 
is sharp. 

\begin{prop}\label{Cex}
Let $0\le \alpha<1$. There exists an operator $P$ on some $L^2(m)$ and $f\in L^2(m)$ such that, for every $\varepsilon>0$,  $\sum_{n\in \N} (\log \log (n+3))^{2-\varepsilon}(n+1)^{1-2\alpha}\|P^nf\|_2^2
<\infty$ and the series $\sum_{n\in \N} (n+1)^{-\alpha}P^n(f)$ diverges $m$-almost everywhere.
\end{prop}
\noindent {\bf Remarks.} The proof is related to some arguments of 
Gaposhkin \cite{Gaposhkin} and makes use of a counterexample by Tandori in the theory 
of orthogonal series. The construction of the operator $P$  is related to the construction of the operator used in the proof of Proposition \ref{prop-eq}. 
Actually, the operator $P$ used in the proof is a one-sided shift, hence is a 
co-isometry which prevent it from being Ritt. This raises the question 
whether it is possible to find a Ritt contraction satisfying the conclusion 
of the Proposition.

\noindent {\bf Proof.} Let $(\varepsilon_n)_{n\in \N}$ be an 
orthormal system on some $L^2(m)$ that we shall specify later. 
We define an operator $P$ on $\overline{{\rm Vect}\{\varepsilon_n\, 
:\, n\in \N\}}$ as follows. For every $f=\sum_{n\in \N}
c_n\varepsilon_n $ let $Pf:= \sum_{n\in \N}c_{n+1}\varepsilon_n$. 
One may extend $P$ to the whole $L^2(m)$ as one please. 

\smallskip

For every $n\in \N$, let $c_n:=\frac1{(n+1)^{3/2-\alpha}\sqrt{\log (n+2)
\log\log(n+3)^{3/2}}}$ and define $f$ as above. 

We have 
$$
\sum_{n=0}^{2^N} (n+1)^{-\alpha}P^n f= \sum_{n=0}^{2^N} 
\sum_{k\ge 0} (n+1)^{-\alpha}c_{k+n}\varepsilon_k
$$
We first prove that 
\begin{equation}\label{first-ser}
\sum_{n=0}^{2^N} 
\sum_{k\ge 2^N+1} (n+1)^{-\alpha}c_{k+n}\varepsilon_k
:= v_{N}\underset{N\to \infty}\longrightarrow 
0 \qquad \mbox{$m$-a.e.}
\end{equation}
We have 
\begin{gather*}
\|v_{N}\|_{L^2(m)}^2= \sum_{k\ge 2^N+1}(\sum_{n=0}^{2^N}(n+1)^{-\alpha}c_{k+n}
)^2\le C 2^{2N(1-\alpha)}\sum_{k\ge 2^N+1}c_{k}^2 
\le \frac{C'}{N(\log (N+1))^3}\, . 
\end{gather*}
Hence, $\sum_{N\in \N}\|v_N\|_{L^2(m)}^2<\infty$ and \eqref{first-ser} 
holds.

Next, we prove that 
\begin{equation}\label{second-ser}
\sum_{k=0}^{2^N} 
\sum_{n\ge 2^N+1} (n+1)^{-\alpha}c_{k+n}\varepsilon_k
:= w_{N}\underset{N\to \infty}\longrightarrow 
0 \qquad \mbox{$m$-a.e.}
\end{equation}

We have 
\begin{gather*}
\|v_{N}\|_{L^2(m)}^2= \sum_{k=0}^{2^N}(\sum_{n\ge 2^N+1}(n+1)^{-\alpha}c_{k+n}
)^2\le C2^{N} (\sum_{n\ge 2^N+1} (n+1)^{-\alpha}c_{n})^2 
\le \frac{C'}{N(\log (N+1))^3}\, . 
\end{gather*}
Hence, $\sum_{N\in \N}\|w_N\|_{L^2(m)}^2<\infty$ and \eqref{second-ser} 
holds.
\smallskip

Combining those first results, we see that we are back to finding 
$(\varepsilon_n)_{n\in \N}$ such that 
$(\sum_{k=0}^{2^N}\sum_{n\ge 0} (n+1)^{-\alpha}c_{n+k}\varepsilon_k)_{N\in \N}
$ diverges $m$-almost everywhere.

\smallskip

For every $k\ge 1$, define 
\begin{gather*}
\alpha_k:= \int_0^\infty
\frac{dx}{x^{\alpha}(x+k)^{3/2-\alpha}
(\log (x+k)(\log \log (x+k+2))^3)^{1/2}
}\\=k^{-1/2}\int_0^\infty \frac{du}{u^{\alpha}(u+1)^{3/2-\alpha}(\log (ku+k)(\log \log (ku+k+2))^3)^{1/2}
}\, .
\end{gather*}

Hence, $(\sqrt k \alpha_k)_{k\ge 1}$ is non-increasing. Moreover, it 
is not hard to see that the series 
$$
\sum_{k\ge 0}\Big(\sum_{n\ge 0} (n+1)^{-\alpha}c_{n+k}-\alpha_k) \varepsilon_k
$$
converges $m$-almost everywhere. 

\smallskip

Then, by a result of Tandori, see Theorem 2.9.1 page 143 of Alexits 
\cite{Alexits} (combined with Theorem 2.7.3 page 120) there exists an orthonormal system 
$(\varepsilon_n)_{n\in \N}$ such that the  $(\sum_{n=0}^{2^N} 
\alpha_n\varepsilon_n)_{N\in \N}$ diverges $m$-a.e., and the proof is 
complete. \hfill $\square$



\section{\bf Extensions, problems}

Recall that an operator $T$ on $\H$ is said to be similar to a contraction 
if there exists a continuous invertible operator $V$ from 
$\H$ onto $\H$ such that $\|VTV^{-1}\|\le 1$, i.e. such that 
$VTV^{-1}$ be a contraction. 

Clearly, all the results from section \ref{section-norm} extend to 
operators that are similar to a contraction. Now, when $\H=L^2(m)$, 
it can be checked that all the results from section \ref{section-ae} 
also hold for operators that are similar to a contraction, even though the 
operator $V$ in the definition need not be positive. 

\medskip

The most general class of operators on $\H$ to which one may hope 
to extend Proposition \ref{prop-alpha}  is the class 
of power bounded operators. Recall that an operator $P$ on 
$\H$ is said to be power bounded if $\sup_{n\in \N}\|P^nf\|<\infty$. 
However, we shall see that this extension is not possible, even if we ask the operator to be Ritt.   
 \medskip

 The next proposition is a reformulation of Proposition 8.2 of \cite{Lemerdy}.

\begin{prop}\label{prop-alpha1}
There exists a Ritt power bounded operator $T$ on some Hilbert space $\H$ such that, taking $\alpha=0$, $(ii)$ of Proposition \ref{prop-alpha} 
does not imply $(i)$.
\end{prop}
\noindent {\bf Proof.} Let $T$ be the operator  defined in Proposition 8.2 of 
Le Merdy \cite{Lemerdy}. Then, $T$ is power bounded and Ritt and has no fixed point. Assume that for every $f\in \H$, the condition 
$\sum_{n\in \N} n\|T^nf\|_\H^2<\infty$ implies that $\sum_{n\in \N}T^nf$ converges in $\H$. Then, arguing as in the proof
of Proposition \ref{prop-eq}, there exists $C>0$ such that 
$\|\sum_{n\in \N}T^nf\|_\H\le C \sum_{n\in \N} n \|T^nf\|_\H^2$, whenever 
the right-hand side converges.

\smallskip

It follows from the proof of Proposition 8.2 of \cite{Lemerdy} that for every $g\in \H$,  $\sum_{n\in \N} n\|(I-T)T^ng\|_\H^2<\infty$. Hence, 
for every $g\in \H$, the series $\sum_{n\in \N}T^n (I-T)g $ converges in $\H$, say to $h$. Then, $\|h\|\le C \sum_{n\in \N} n\|(I-T)T^ng\|_\H^2$ and 
$(I-T)h=(I-T)g$. Since $T$ has no fixed point, we infer that $g=h$ and 
$\|g\|\le C \sum_{n\in \N} n\|(I-T)T^ng\|_\H^2$ for every $g\in \H$. 
But it is proved in \cite{Lemerdy} that this cannot hold.  \hfill $\square$

\medskip

We now give an extension of Corollary \ref{cor} to the case 
where $\alpha=1$.

\begin{prop}
Let $P$ be a contraction on $L^2(m)$. Let $f\in L^2(m)$ be such that 
\begin{equation}\label{EHT}
\sum_{n\in \N} \frac{\log( n+1)}{n+1}\|P^nf\|_2^2<\infty\, .
\end{equation}
Then, $\sum_{n \in \N} \frac{P^nf}{n+1}$ converges in $L^2(m)$. 
If moreover 
$$
\sum_{n\in \N} \frac{\log( n+1)\big(\log \log \log (n+9)\big)^2}{n+1}\|P^nf\|_2^2<\infty\, 
$$ 
then, $\sum_{n \in \N} \frac{P^nf}{n+1}$ converges $m$-almost everywhere 
and the associated maximal function is in $L^2(m)$.
\end{prop}
\noindent {\bf Proof.} The norm convergence follows easily from 
Proposition \ref{prop-norm}. We now give the main argument for 
the proof of the almost everywhere convergence. Then, the rest of 
the proof is similar to that of Proposition \ref{theo-ae}. Assume 
\eqref{EHT}.
We have
$$
\|\sum_{k=2^{2^n}}^{2^{2^{n+1}-1}} \frac{|P^kf|}{k}\|_{L^2(m)}^2\le 
C2^{2n}\|P^{2^{2^n}}f\|_{L^2(m)}^2\, .
$$
It is not hard to see that \eqref{EHT} implies that $\sum_{n\ge 0}2^{2n} \|P^{2^{2^n}}\|_{L^2(m)}^2<\infty$. 
Hence it suffices to prove the desired 
convergence of the series along the sequence $(2^{2^n})_{n\in \N}$. 
\hfill $\square$

\medskip 

A natural question is the following : does there 
exist analogous results to, say, Proposition \ref{prop-alpha} 
for contractions of $L^p$ spaces. For instance,  by Cohen, Cuny and Lin \cite{CCL}, if $P$ is a positive Ritt contraction 
of some $L^p$, $1<p\le 2$ then the condition $\sum_{n\in \N} 
(n+1)^{1-p\alpha }\|P^nf\|_p^p<\infty$ is sufficient for the convergence 
in $L^p$ of $\sum_{n\in \N} (n+1)^{-\alpha}P^nf$ (and the a.e. convergence holds as well). The approach used in the present paper partially work 
for Markov operators. However, it does not seem to  allow one to extend 
the results of \cite{CCL} to Markov operators. It would be interesting 
either to prove that $L^p$ extensions are possible or 
to find an example where it cannot.

\begin{appendix}

\section{}
We made use of the following lemma which is related to Lemma 2.7 of Peligrad and Utev \cite{PU}.

\begin{lem}\label{subadditive}
Let $(V_n)_{n\ge 1}$ be a nonnegative subadditive sequence 
(i.e. $V_{n+m}\le V_n+V_m$  for every $m,n\ge 1$). Then, for every 
$q\ge 1$ and every $p>1$, there exists $C>0$ such that    
$\sum_{n\ge 1}\frac{\max_{1\le i\le n}V_n^q}{n^p}
\le C \sum_{n\ge 0}\frac{V_{2^n}^q}{2^{np}}$.
\end{lem}
\noindent {\bf Proof. } The proof basically follows the arguments to prove Lemma 4.1 of \cite{PU}. We start with the following basic (dyadic) 
decomposition. For every $r\ge 0$, and every $2^r\le n\le 2^{r+1}-1$, 
we have, by an easy induction, 
$$
V_n \le V_{n-2^r}+V_{2^r}\le \sum_{k=0}^r V_{2^k}\, .
$$
Hence,

\begin{gather*}
\sum_{n\ge 1}\frac{\max_{1\le i\le n}V_n^q}{n^p}
\le \sum_{r\ge 0} \sum_{n=2^r}^{2^{r+1}-1}\frac{\max_{1\le i\le n}V_n^q}{n^p} 
\le \sum_{r\ge 0}\frac1{2^{(p-1)r}}(\sum_{k=0}^rV_{2^k})^q
\end{gather*}

When $q=1$ the result follows. Assume that $q>1$. Let $0<\varepsilon <p-1$.   Using H\"older's inequality (with $1/q+1/q'=1$)
1, we have
\begin{gather*}
\sum_{n\ge 1}\frac{\max_{1\le i\le n}V_n^q}{n^p}\le \sum_{r\ge 0}\frac1{2^{(p-1)r}}(\sum_{j=0}^r2^{k\varepsilon q'q})^{q/q'}(\sum_{k=0}^r 
2^{-k\varepsilon}V_{2^k}^q)
\le C  \sum_{n\ge 0}\frac{V_{2^n}^q}{2^{np}}\, .
\end{gather*}

\end{appendix}

\end{document}